\documentclass[11pt]{amsart}
\usepackage{amssymb}
\usepackage{amsmath}
\usepackage{graphicx,psfrag}

\title{Embedding large subgraphs into dense graphs}

\date{}
\author{Daniela K\"uhn and Deryk Osthus}

\newtheorem{firstthm}{Proposition}
\newtheorem{theorem}[firstthm]{Theorem}
\newtheorem{prop}[firstthm]{Proposition}
\newtheorem{lemma}[firstthm]{Lemma}
\newtheorem{cor}[firstthm]{Corollary}

\newtheorem{conj}[firstthm]{Conjecture}

\def\noproof{{\unskip\nobreak\hfill\penalty50\hskip2em\hbox{}\nobreak\hfill%
       $\square$\parfillskip=0pt\finalhyphendemerits=0\par}\goodbreak}
\def\endproof{\noproof\bigskip}
\newdimen\margin   % needed for macros \textdisplay & \ltextdisplay
\def\textno#1&#2\par{%
   \margin=\hsize
   \advance\margin by -4\parindent
          \setbox1=\hbox{\sl#1}%
   \ifdim\wd1 < \margin
      $$\box1\eqno#2$$%
   \else
      \bigbreak
      \hbox to \hsize{\indent$\vcenter{\advance\hsize by -3\parindent
      \it\noindent#1}\hfil#2$}%
      \bigbreak
   \fi}
\def\proof{\removelastskip\penalty55\medskip\noindent{\bf Proof. }}

\begin{document}
\maketitle

\def\COMMENT#1{}
\def\TASK#1{}

\def\eps{{\varepsilon}}
\newcommand{\ex}{\mathbb{E}}
\newcommand{\pr}{\mathbb{P}}
\newcommand{\hcf}{{\rm hcf}}
\newcommand{\D}{\mathcal D}

\begin{abstract} \noindent
What conditions ensure that a graph $G$ contains some given spanning subgraph $H$?
The most famous examples of results of this kind are probably Dirac's theorem on Hamilton cycles and 
Tutte's theorem on perfect matchings.
Perfect matchings are generalized by perfect $F$-packings, where instead of covering 
all the vertices of $G$ by disjoint edges, we want to cover $G$ by disjoint copies of a (small) graph $F$.
It is unlikely that there is a characterization of all graphs $G$ which contain a perfect $F$-packing,
so as in the case of Dirac's theorem it makes sense to study conditions on the minimum degree
of~$G$ which guarantee a perfect $F$-packing. 

The Regularity lemma of Szemer\'edi and the Blow-up lemma of Koml\'os, S\'ark\"ozy and Szemer\'edi 
have proved 
to be powerful tools in attacking such problems and quite recently, several long-standing problems 
and conjectures in the area have been solved using these. 
In this survey, we give an outline of recent progress 
(with our main emphasis on $F$-packings, Hamiltonicity problems and tree embeddings)
and describe some of the methods involved.
\end{abstract}

\section{Introduction, overview and basic notation}\label{intro}
In this survey, we study the question of when a graph $G$ contains some given large or spanning
graph $H$ as a subgraph. Many important problems can be phrased in this way:
one example is Dirac's theorem, which states that every graph $G$ on $n \ge 3$ vertices with
minimum degree at least $n/2$ contains a Hamilton cycle. 
Another example is Tutte's theorem on perfect matchings which gives a characterization of all those graphs
which contain a perfect matching (so $H$ corresponds to a perfect matching in this case).
A result which gives a complete characterization of all those graphs $G$ which contain $H$
(as in the case of Tutte's theorem) is of course much more desirable than a sufficient condition
(as in the case of Dirac's theorem). However, for most $H$ that we consider, it is unlikely that
such a characterization exists as the corresponding decision problems are usually NP-complete.
So it is natural to seek simple sufficient conditions.
Here we will focus mostly on degree conditions. This means that $G$ will usually be a dense graph and that
we have to restrict $H$ to be rather sparse in order to get interesting results.
We will survey the following topics: 
\begin{itemize}
\item a generalization of the matching problem, which is called the \emph{$F$-packing} 
or \emph{$F$-tiling} problem (here the aim is to cover the vertices of $G$ with 
disjoint copies of a fixed graph $F$ instead of disjoint edges);
\item Hamilton cycles (and generalizations) in graphs, directed graphs and hypergraphs;
\item large subtrees of graphs;
\item arbitrary subgraphs $H$ of bounded degree;
\item Ramsey numbers of sparse graphs.
\end{itemize}
A large part of the progress in the above areas is due to the Regularity lemma of Szemer\'edi~\cite{reglem} and
the Blow-up lemma of Koml\'os, S\'ark\"ozy and Szemer\'edi~\cite{KSSblowup}. Roughly speaking, the former states that one can decompose an
arbitrary large dense graph into a bounded number of random-like graphs. The latter is a powerful tool for embedding
spanning subgraphs $H$ into such random-like graphs. In the final section we give a formal statement of these
results and describe in detail an application to a special case of the $F$-packing problem.
We hope that readers who are unfamiliar with these tools will find this a useful guide to how they can be applied.

There are related surveys in the area by Koml\'os and Simonovits~\cite{KSi} 
(some minor updates were added later in~\cite{KSSS})
and by Koml\'os~\cite{JKblowup}.
However, much has happened since these were written and the emphasis is different in each case.
So we hope that the current survey will be a useful complement and update to these.
In particular, as the title indicates, our focus is mainly on embedding large subgraphs
and we will ignore other aspects of regularity/quasi-randomness.
There is also a recent survey on $F$-packings (and so-called $F$-decompositions) by Yuster~\cite{yuster},
which is written from a computational perspective.

%%%%%%%%%%%%%%%%%%%%%%%%%%%%%%%%%%%%%%%%%%%%%%%%%%%%%%%%%%%%%%%%%%%%%%%%%%%%%%%%%%%%%%%%%%%%
\section{Packing small subgraphs in graphs} \label{packings}

\subsection{$F$-packings in graphs of large minimum degree}

Given two graphs $F$ and
$G$, an \emph{$F$-packing in $G$} is a collection of vertex-disjoint copies
of $F$ in $G$. (Alternatively, this is often called an $F$-tiling.)
$F$-packings are natural generalizations of graph matchings
(which correspond to the case when $F$ consists of a single edge). 
An $F$-packing in $G$ is called \emph{perfect} if it
covers all vertices of $G$. In this case, we also say that $G$ contains
an \emph{$F$-factor} or a \emph{perfect $F$-matching}.
If $F$ has a component which contains at least~3 vertices then the question whether 
$G$ has a perfect $F$-packing is difficult from both a structural and algorithmic
point of view: Tutte's theorem characterizes those graphs which have a perfect 
$F$-packing if $F$ is an edge but for other connected graphs~$F$ no such characterization
is known. Moreover, Hell and Kirkpatrick~\cite{HKsiam}
showed that the decision problem of whether a graph $G$ has a perfect $F$-packing
is NP-complete if and only if $F$ has a component which contains at least~3 vertices.
So as mentioned earlier, this means that it makes sense to search for degree
conditions which ensure the existence of a perfect $F$-packing.
The fundamental result in the area is the Hajnal-Szemer\'edi theorem:
\begin{theorem}{\bf (Hajnal and Szemer\'edi~\cite{HSz})} \label{hajnalsz}
Every graph whose order~$n$ is divisible by~$r$ and whose minimum degree is at least
$(1-1/r)n$ contains a perfect $K_r$-packing. 
\end{theorem}
The minimum degree condition is easily seen to be best possible. (The case when $r=3$
was proved earlier by Corr\'adi and Hajnal~\cite{CH}.) 
The result is often phrased in terms of colourings: any graph~$G$ whose order is
divisible by~$k$ and with $\Delta(G) \le k-1$
has an equitable $k$-colouring, i.e.~a colouring with colour classes of equal size.
(So $k:=n/r$ here.)
Theorem~\ref{hajnalsz} raises the question of what minimum degree condition forces a perfect 
$F$-packing for arbitrary graphs~$F$. The following result gives a general bound.
\begin{theorem}\label{KSS}{\bf (Koml\'os, S\'ark\"ozy and Szemer\'edi~\cite{KSSz01})}
For every graph $F$ 
there exists a constant $C=C(F)$ such that every
graph $G$ whose order $n$ is divisible by $|F|$ and whose
minimum degree is at least $(1-1/\chi(F))n+C$ contains a perfect $F$-packing.
\end{theorem}
This confirmed a conjecture of Alon and Yuster~\cite{AY96}, who had obtained the
above result with an additional error term of~$\eps n$ in the minimum degree condition.
As observed in~\cite{AY96}, there are graphs $F$ 
for which the above constant~$C$ cannot be omitted completely (e.g.~$F=K_{s,s}$ where $s\ge 3$
and $s$ is odd).
Thus one might think that this settles the question of which minimum degree 
guarantees a perfect $F$-packing. However, we shall see that this is \emph{not} the case.
There are graphs $F$ for which the bound on the minimum 
degree can be improved significantly: we can often replace 
$\chi(F)$ by a smaller parameter. For a detailed statement of this, 
we define the \emph{critical chromatic number} $\chi_{cr}(F)$ of a graph $F$ as 
$$
\chi_{cr}(F):=(\chi(F)-1)\frac{|F|}{|F|-\sigma(F)},$$ 
where $\sigma(F)$ denotes
the minimum size of the smallest colour class in an optimal colouring
of $F$. (We say that a colouring of $F$ is \emph{optimal} if it uses exactly $\chi(F)$ colours.)
So for instance a $k$-cycle $C_k$ with $k$ odd has
$\chi_{cr}(C_k)=2+2/(k-1)$.
Note that $\chi_{cr}(F)$ always satisfies  $\chi(F)-1 < \chi_{cr}(F) \le \chi(F)$  
and equals $\chi(F)$ if and only if for every optimal colouring of $F$ 
all the colour classes have equal size.
The critical chromatic number was introduced by Koml\'os~\cite{JKtiling}.
He (and independently Alon and Fischer~\cite{AF99}) observed that for \emph{any}
graph~$F$ it can be used to give a lower bound on the minimum degree that guarantees
a perfect $F$-packing.%

\begin{prop}\label{propKomlos}
For every graph $F$ and every integer $n$ that is divisible by $|F|$ there exists
a graph $G$ of order $n$ and minimum degree $\lceil(1-1/\chi_{cr}(F))n\rceil-1$
which does not contain a perfect $F$-packing.
\end{prop} 
Given a graph~$F$, the graph $G$ in the proposition is constructed as follows:
write $k:=\chi(F)$ and let $\ell \in \mathbb{N}$ be arbitrary. 
$G$ is a complete $k$-partite graph with vertex classes $V_1,\dots,V_k$,
where $|V_1|=\sigma(F)\ell-1$, $n=\ell |F|$ and the sizes of $V_2,\dots,V_k$ are as equal
as possible. Then any perfect $F$-packing would consist of $\ell$ copies of $F$.
On the other hand, each such copy would contain at least $\sigma(F)$ vertices in $V_1$, which 
is impossible.

Koml\'os also showed that the critical chromatic number is the parameter which governs the 
existence of \emph{almost} perfect packings in graphs of large minimum degree.
(More generally, he also determined the minimum degree which ensures that a given fraction
of vertices is covered.)
\begin{theorem}\label{thmKomlos}{\bf (Koml\'os~\cite{JKtiling})}
For every graph $F$ and every $\gamma>0$ there exists an integer
$n_0=n_0(\gamma,F)$ such that every graph $G$ of order $n\ge n_0$
and minimum degree at least $(1-1/\chi_{cr}(F))n$ contains an $F$-packing
which covers all but at most $\gamma n$ vertices of~$G$.
\end{theorem}
By making $V_1$ slightly smaller in the previous example, it is easy to see that the 
minimum degree bound in Theorem~\ref{thmKomlos} is also best possible.
Confirming a conjecture of Koml\'os~\cite{JKtiling}, Shokoufandeh and Zhao~\cite{SZ,SZ3} subsequently
proved that the number of uncovered vertices can be reduced to a constant depending only on~$F$.%

We~\cite{KOmatch} proved that for any graph~$F$, either its critical chromatic number or its
chromatic number is the relevant parameter which governs the existence of perfect 
packings in graphs of large minimum degree. The classification depends on a
parameter which we call the \emph{highest common factor} of $F$.

This is defined as follows for non-bipartite graphs~$F$.
Given an optimal colouring $c$ of~$F$, let $x_1\le x_2\le \dots\le x_{\ell}$
denote the sizes of the colour classes of~$c$. Put
$\D (c):= \{ x_{i+1}-x_i\,|\, i=1,\dots, \ell-1 \}.$
Let $\D(F)$ denote the union of all the sets $\D (c)$ taken over all optimal colourings $c$.
We denote by $\hcf(F)$ the highest common factor of all integers in $\D(F)$.
If $\D(F)=\{0\}$ we set $\hcf(F):=\infty$.
Note that if all the optimal colourings of $F$ have the property that
all colour classes have equal
size, then $\D(F)=\{0\}$ and so $\hcf(F) \neq 1$ in this case.
In particular, if $\chi_{cr}(F)=\chi(F)$, then $\hcf(F) \neq 1$.
So for example, odd cycles of length at least 5 have $\hcf=1$ whereas complete graphs 
have $\hcf \neq 1$.

The definition can be extended to bipartite graphs $F$.
For connected bipartite graphs, we always have $\hcf(F) \neq 1$, but for disconnected bipartite graphs
the definition also takes into account the relative sizes of the components of~$F$ (see~\cite{KOmatch}).

We proved that in Theorem~\ref{KSS}
one can replace the chromatic number by the critical chromatic number if $\hcf(F)=1$. 
(A much simpler proof of a weaker result can be found in~\cite{KOSODA}.)

\begin{theorem}{\bf (K\"uhn and Osthus~\cite{KOmatch})} \label{thmmain}
Suppose that $F$ is a graph with
$\hcf(F)=1$. Then there exists
a constant $C=C(F)$ such that every graph $G$ whose order $n$ is divisible by~$|F|$
and whose minimum degree is at least $(1-1/\chi_{cr}(F))n+C$
contains a perfect $F$-packing.
\end{theorem}

Note that Proposition~\ref{propKomlos} shows that the result is best possible up to the
value of the constant~$C$. A simple modification of the examples in~\cite{AF99,JKtiling}
shows that there are graphs~$F$ for which the constant $C$ cannot be omitted
entirely.
Moreover, it turns out that Theorem~\ref{KSS} is already best possible up to
the value of the constant~$C$ if $\hcf(F)\neq 1$. To see this, for simplicity assume that
$k:=\chi(F) \ge 3$ and $n=k \ell |F|$ for some $\ell \in \mathbb{N}$
and let $G$ be a complete $k$-partite graph with vertex classes $V_1,\dots,V_k$, where
$|V_1|:=\ell |F|-1$, $|V_2|:=\ell |F|+1$ and $|V_i|=\ell |F|$ for $i \ge 3$. Consider any $F$-packing 
$F_1,\dots,F_t$ in~$G$. Let $G_i$ be the graph obtained from $G$ by removing $F_1,\dots, F_{i}$.
So $G=G_0$. If $t:=\hcf(F) \neq 1$, then the vertex classes $V_i^1$ of  $G_1$ still have property that
$|V_1^1|-|V_k^1| \not\equiv 0$ modulo $t$. More generally, this property is preserved for all
$G_i$, so the original $F$-packing cannot cover all the vertices in $V_1 \cup V_k$.

One can now  combine Theorems~\ref{KSS} and~\ref{thmmain} 
(and the corresponding lower bounds which are discussed
in detail in~\cite{KOmatch}) to obtain a complete answer to the
question of which minimum degree forces a perfect $F$-packing (up to an additive constant). 
For this, let 
$$\chi^*(F):=
\begin{cases}
\chi_{cr}(F) &\text{ if $\hcf (F)=1$};\\
\chi(F) &\text{ otherwise}.
\end{cases}
$$
Also let $\delta(F,n)$ denote the smallest integer $k$ such that every graph $G$ 
whose order $n$ is divisible by~$|F|$ and with $\delta(G)\ge k$ contains a perfect $F$-packing.   
\begin{theorem}{\bf (K\"uhn and Osthus~\cite{KOmatch})}\label{thmmaingeneral}
For every graph $F$ there exists a constant $C=C(F)$ such that
$$\left( 1-\frac{1}{\chi^*(F)} \right)n-1\le \delta(F,n)
\le \left(1-\frac{1}{\chi^*(F)} \right)n+C.$$
\end{theorem}
The constant~$C$ appearing in Theorems~\ref{thmmain} and~\ref{thmmaingeneral} is rather large
since it is related to the number of partition classes (clusters) obtained by the Regularity lemma.
It would be interesting to know whether one can take e.g.~$C=|F|$ (this holds
for large~$n$ in Theorem~\ref{KSS}). Another open problem
is to characterize all those graphs~$F$ for which $\delta(F,n)=\lceil(1-1/\chi^*(F))n\rceil$.
This is known to be
the case for complete graphs by Theorem~\ref{hajnalsz} and all graphs with at most
$4$ vertices (see~Kawarabayashi~\cite{ken} for a proof of the case when $F$ is a $K_4$ minus an edge 
and a discussion of the other cases). If $n$ is large, this is also known to hold
for cycles (this follows from Theorem~\ref{abbasi} below) and
for the case when $F$ is a complete graph minus an edge~\cite{KOKlminus} (the latter was conjectured in~\cite{ken}). %

\subsection{Ore-type degree conditions}
Recently, a simple proof (based on an inductive argument)
of the Hajnal-Szemer\'edi theorem was found by Kierstead and 
Kostochka~\cite{KK1}.
Using similar methods, they subsequently strengthened this to an Ore-type condition~\cite{KK2}:
\begin{theorem}{\bf (Kierstead and Kostochka~\cite{KK2})}
Let $G$ be a graph whose order~$n$ is divisible by~$r$. If $d(x)+d(y) \ge 2(1-1/r) n-1$ for
all pairs $x\neq y$ of nonadjacent vertices, then~$G$ has a perfect~$K_r$-packing.
\end{theorem}
Equivalently, if a graph~$G$ whose order is divisible by~$k$ satisfies $d(x)+d(y) \le 2k-1$
for every edge~$xy$, then~$G$ has an equitable~$k$-colouring. (So $k:=n/r$.)
Recently, together with Treglown~\cite{KOTore}, we proved an Ore-type analogue of Theorem~\ref{thmmaingeneral}
(but with a linear error term $\varepsilon n$ instead of the additive constant~$C$).
The result in this case turns out to be genuinely different: 
again, there are some graphs $F$ for which the degree condition depends on $\chi(F)$ and some
for which it depends on $\chi_{cr}(F)$. However, there are also graphs~$F$ for which it depends on a 
parameter which lies strictly between $\chi_{cr}(F)$ and $\chi(F)$.
This parameter in turn depends on how many additional colours are necessary to extend colourings
of neighbourhoods of certain vertices of~$F$ to a colouring of~$F$.
It is an open question whether the linear error term in~\cite{KOTore} can be reduced to a constant one.

\subsection{$r$-partite versions}
Also, it is natural to consider $r$-partite versions of the Hajnal-Szemer\'edi theorem.
For this, given an~$r$-partite graph~$G$, let $\delta'(G)$ denote the minimum over all vertex
classes~$W$ of~$G$ and all vertices~$x\notin W$ of the number of neighbours of~$x$ in~$W$.
The obvious question is what value of $\delta'(G)$ ensures that~$G$ has a perfect~$K_r$-packing.
The following (surprisingly difficult) conjecture is implicit in~\cite{MM}.
Fischer~\cite{Fisher} originally made a stronger conjecture which did not include the
`exceptional' graph~$\Gamma_{r,n}$ defined below.
\begin{conj} \label{partite}
Suppose that $r\ge 2$ and that~$G$ is an $r$-partite graph with vertex classes of size $n$.
If $\delta'(G) \ge (1-1/r)n$, then $G$ has a perfect $K_r$-packing
unless both $r$ and $n$ are odd and $G=\Gamma_{r,n}$.
\end{conj}
To define the graph $\Gamma_{r,n}$, we first construct a graph $\Gamma_r$: 
its vertices are labelled $g_{ij}$ with $1 \le i,j \le r$.
We have an edge between $g_{ij}$ and $g_{i'j'}$
if $i \neq i'$, $j \neq j'$ and $j \le r-2$ or $j' \le r-2$.
We also have an edge if $i \neq i'$ and we have either $j=j'=r-1$ or $j=j'=r$ (see Fig.~1).
\begin{figure}\label{g33}
\centering
\includegraphics[scale=0.32]{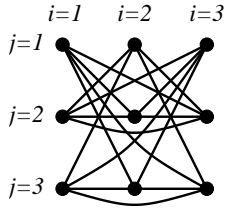}
\caption{The graph $\Gamma_3=\Gamma_{3,1}$ in Conjecture~\ref{partite}}
\end{figure}
$\Gamma_{r,n}$ is then obtained from $\Gamma_r$ by replacing each vertex with an independent set
of size $n/r$ and replacing each edge with a complete bipartite graph.

To see that $\Gamma_{r,n}$ has no perfect $K_r$-packing when both $r$ and $n$ are odd, 
let $W_\ell$ denote the set of vertices of $\Gamma_{r,n}$ which correspond to a vertex of
$\Gamma_{r}$ with $j =\ell$. Note that every copy of~$K_r$ which covers a vertex in
$W_1\cup\dots\cup W_{r-2}$ has to contain at least 2 vertices in~$W_{r-1}$ or
at least~2 vertices in~$W_r$. So in order to cover all vertices in $W_1\cup\dots\cup W_{r-2}$
we can only use copies of~$K_r$ which contain exactly~2 vertices in~$W_{r-1}$ or
exactly~2 vertices in~$W_r$. But since $|W_{r-1}|=|W_r|=n$ is odd this means that it is impossible to
cover all vertices of $\Gamma_{r,n}$ with vertex-disjoint copies of~$K_r$.
(Note that the argument uses only that~$n$ is odd, but we cannot have that~$n$ is odd and~$r$ is even.)

A much simpler example which works for all~$r$ and~$n$ but which gives a weaker bound when $r$ and $n$ are odd 
is obtained as follows: choose a set $A$ which has
less than $(1-1/r)n$ vertices in each vertex class and include all edges which have at least one endpoint in~$A$.
For large $n$, the case $r=3$ of Conjecture~\ref{partite} was solved by Magyar and Martin~\cite{MM} and
the case $r=4$ by Martin and Szemer\'edi~\cite{MS}, both using the Regularity lemma
(the case $r=2$ is elementary). Johansson~\cite{johansson} had earlier proved an approximate version of the case~$r=3$.
Csaba and Mydlarz~\cite{csabamulti} proved a result which implies that Conjecture~\ref{partite} holds 
approximately when~$r$ is large (and $n$ large compared to~$r$). 
Generalizations to packings of arbitrary graphs were considered in~\cite{hladkyschacht,zhaomartin,zhaobip}.
A variant of the problem (where one considers usual minimum degree $\delta(G)$) was considered by
Johansson, Johansson and Markstr\"om~\cite{JJM}. They solved the case $r=3$ and gave bounds for the case
$r>3$. This problem is related to bounding the so-called `strong chromatic number'.

\subsection{Hypergraphs} \label{hyperpack}
(Perfect) $F$-packings have also been investigated for the case when~$F$
is a uniform hypergraph. Unsurprisingly, the hypergraph problem turns out to be much more difficult than the graph problem. There are two natural notions of (minimum) degree
of the `dense' hypergraph~$G$. Firstly, one can consider the vertex degree.
Secondly, given an $r$-uniform hypergraph~$G$ and an $(r-1)$-tuple~$W$ of vertices in $G$, the 
degree of $W$ is defined
to be the number of hyperedges which contain~$W$.
This notion of degree is called  \emph{collective degree} or 
\emph{co-degree}. 
In contrast to the graph case, even the minimum collective degree which ensures a perfect matching (i.e.~when~$F$
consists of a single edge) is not easy to determine. 
R\"odl, Ruci\'nski and Szemer\'edi~\cite{RRS} gave a precise solution to this problem, the answer
turns out to be close to $n/2$.
This improved bounds of~\cite{KOhypermatch,RRSapprox}.
An $r$-partite version (which is best possible for infinitely many values of~$n$)
was proved by Aharoni, Georgakopoulos and Spr\"ussel~\cite{AGS}.
The minimum vertex degree which forces the existence of a perfect matching is unknown.
It is natural to make the following conjecture (a related $r$-partite version is conjectured in~\cite{AGS}).%
\begin{conj} \label{hypermatch}
For all integers $r$ and all $\eps>0$ there is an integer $n_0=n_0(r,\eps)$ so that the following holds for all
$n\ge n_0$ which are divisible by~$r$: if~$G$ 
is an $r$-uniform hypergraph on $n$ vertices whose minimum vertex degree is at least 
$$(1-(1-1/r)^{r-1}+\eps)\binom{n}{r-1},$$ then $G$ has a perfect matching.
\end{conj}
The following construction gives a corresponding lower bound: let $V$ be a set of $n$ vertices
and let $A \subseteq V$ be a set of less than $n/r$ vertices and include as  
hyperedges all $r$-tuples with at least one vertex in $A$.
The case $r=3$ of the conjecture was proved recently by Han, Person and Schacht~\cite{HPS}.

A hypergraph analogue of Theorem~\ref{thmmaingeneral} currently seems out of reach.
So far, the only hypergraph~$F$ (apart from the single edge) for which the approximate minimum collective degree
which forces a perfect $F$-packing has been determined is the $3$-uniform hypergraph with 4 vertices and
2 edges~\cite{KOloose}.
Pikhurko~\cite{pikhurko} gave bounds on the minimum collective degree which forces the complete $3$-uniform hypergraph
on $4$ vertices. In the same paper, he also shows that if $\ell \ge r/2$ and $G$ is an $r$-uniform hypergraph
where every $\ell$-tuple of vertices is contained
in at least $(1/2+o(1))\binom{n}{r-\ell}$ hyperedges, then $G$ has a perfect matching, which is best possible 
up to the $o(1)$-term. This result is rather surprising in view of the fact that
Conjecture~\ref{hypermatch} (which corresponds to the case when $\ell=1$) has a rather different form.
Further results on this question are also proved in~\cite{HPS}.

%%%%%%%%%%%%%%%%%%%%%%%%%%%%%%%%%%%%%%%%%%%%%%%%%%%%%%%%%55

\section{Trees}

One of the earliest applications of the Blow-up lemma was the solution
by Koml\'os, S\'ark\"ozy and Szemer\'edi~\cite{KSStrees1} 
of a conjecture of Bollob\'as on the existence of given bounded degree spanning trees. 
The authors later relaxed the condition of bounded degree to obtain the following result.
\begin{theorem} {\bf (Koml\'os, S\'ark\"ozy and Szemer\'edi~\cite{KSSztrees})}
For any $ \gamma > 0 $  there exist constants $c>0$ and $n_0$ with the following properties. 
If $ n \geq n_0 $, $T$ is a tree of order $n$ with $ \Delta (T) \leq { cn /\log n}$, 
and $G$ is a graph of order $n$ with $ { \delta (G)} \geq (1/2+ { \gamma })n $, then $T$ is a subgraph of $G$.
\end{theorem}
The condition $ \Delta (T) \leq { cn /\log n}$ is best possible up to the value of $c$.
(The example given in~\cite{KSSztrees} to show this is a random graph $G$ with edge probability $0.9$
and a tree of depth $2$ whose root has degree close to $\log n$.)

It is an easy exercise to see that every graph of minimum degree at least $k$ contains
any tree with $k$ edges. The following classical conjecture would imply that we can replace
the minimum degree condition by one on the average degree.
\begin{conj} {\bf (Erd\H{o}s and S\'os~\cite{ErdosSos})} \label{erdossos}
Every graph of average degree greater than $k-1$ contains any tree with $k$ edges.
\end{conj}
This is trivially true for stars. (On the other hand, stars also show that the bound is best possible in
general.) It is also trivial if one assumes an extra factor of~2
in the average degree. It has been proved for some special classes of trees, most 
notably those of diameter at most~4~\cite{mclennan}. 
%Ajtai, Koml\'os, Simonovits and Szemer\'edi 
%announced a proof for $k$ sufficiently large which is based on the Regularity lemma 
%(i.e.~there is a constant $k_0$ so that the conjecture holds
%for all trees with at least $k_0$ edges).
The conjecture is also true for `locally sparse' graphs
-- see Sudakov and Vondrak~\cite{vondraksudakov} for a discussion of this.

The following result proves (for large $n$) a related conjecture of Loebl. An approximate 
version was proved earlier by Ajtai, Koml\'os and Szemer\'edi~\cite{AKSloebl}.
\begin{theorem}{\bf (Zhao~\cite{zhao})} \label{zhao}
There is an integer $n_0$ so that every graph $G$ on $n \ge n_0$ vertices which has
at least $n/2$ vertices of degree at least $n/2$ contains all trees with at most $n/2$ edges.
\end{theorem}
This would be generalized by the following conjecture. 
\begin{conj} {\bf (Koml\'os and S\'os)} \label{KomlosSos}
Every graph $G$ on $n$ vertices which has
at least $n/2$ vertices of degree at least $k$ contains all trees with $k$ edges.
\end{conj}
Again, the conjecture is trivially true (and best possible) for stars.
Piguet and Stein~\cite{PS} proved an approximate version for the case when $k$ is linear in $n$ and $n$ is large.
Cooley~\cite{cooley} as well as Hladk\'y and Piguet~\cite{hladkypiguet} proved an exact version for this case. 
All of these proofs are based on the Regularity lemma. As with Conjecture~\ref{erdossos}, there
are several results on special cases which are not based on the Regularity lemma.
For instance, Piguet and Stein proved it for trees of diameter at most~5~\cite{PSdiam}.

%%%%%%%%%%%%%%%%%%%%%%%%%%%%%%%%%%%%%%%%%%%%%%%%%%%%%%%%%%%%%%%%%%%%%%%%%%%%%%%%%%%%%%%55

\section{Hamilton cycles}

\subsection{Classical results for graphs and digraphs}
As mentioned in the introduction, the decision problem of whether 
a graph has a Hamilton cycle is NP-complete, so it makes sense to ask for degree
conditions which ensure that a graph has a Hamilton cycle.
One such result is the classical theorem of Dirac.
\begin{theorem}{\bf (Dirac~\cite{dirac})} \label{dirac}
Every graph on $n\ge 3$ vertices with minimum degree
at least $n/2$ contains a Hamilton cycle. 
\end{theorem}
For an analogue in directed graphs it is
natural to consider the \emph{minimum semidegree~$\delta^0(G)$} of a digraph $G$,
which is the minimum of its minimum outdegree~$\delta^+(G)$ and its minimum
indegree~$\delta^-(G)$. (Here a directed graph may have two edges between a pair of vertices, 
but in this case their directions must be opposite.)
The corresponding result is a theorem of Ghouila-Houri~\cite{gh}.
\begin{theorem}{\bf (Ghouila-Houri~\cite{gh})} \label{hamGH}
Every digraph on $n$ vertices with minimum semidegree at least $n/2$ contains a Hamilton cycle. 
\end{theorem}
In fact, Ghouila-Houri proved the stronger result that every strongly connected digraph
of order~$n$ where every vertex has total degree at least~$n$ has a Hamilton cycle.
(When referring to paths and cycles in directed
graphs we always mean that these are directed, without mentioning this explicitly.)
All of the above degree conditions are best possible. 
Theorems~\ref{dirac} and~\ref{hamGH} were generalized to a degree condition on 
pairs of vertices for graphs as well as digraphs:
\begin{theorem}{\bf (Ore~\cite{ore})} 
Suppose that $G$ is a graph with $n \ge 3$ vertices such that every pair $x\neq y$
of nonadjacent vertices satisfies $d(x)+d(y) \ge n$. Then $G$ has a Hamilton cycle.
\end{theorem}
\begin{theorem}{\bf (Woodall~\cite{woodall})} \label{woodall}
Let $G$ be a strongly connected digraph on $n \ge 2$ vertices.
If $d^+(x)+d^-(y)\ge n$ for every pair $x\neq y$ of vertices for which there is no edge from $x$ to $y$,
then $G$ has a Hamilton cycle.
\end{theorem}

There are many generalizations of these results. The survey~\cite{gould} gives an overview
for undirected graphs
and the monograph~\cite{digraphsbook} gives a discussion of directed versions.
Below, we describe some recent progress on degree conditions for Hamilton cycles, much 
of which is based on the Regularity lemma.% 
\COMMENT{mention algos (Sarkozy) somewhere}

\subsection{Hamilton cycles in oriented graphs}
Thomassen~\cite{thomassen_79}
raised the natural question of determining the minimum
semidegree that forces a Hamilton cycle in
an \emph{oriented graph} (i.e.~in a directed graph that can be obtained from a simple
undirected graph by orienting its edges).
Thomassen initially believed that the correct minimum semidegree bound should be $n/3$
(this bound is obtained by considering a `blow-up' of an oriented triangle).
However, H\"aggkvist~\cite{HaggkvistHamilton} later gave the following construction which 
gives a lower bound of $\lceil (3n-4)/8 \rceil -1$ (see Fig.~2).
\begin{figure}\label{extremal2}
\centering\footnotesize
\includegraphics[scale=0.45]{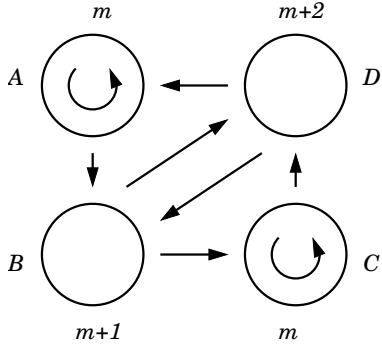}
\caption{An extremal example for Theorem~\ref{main}}
\end{figure}
For $n$ of the form $n=4m+3$ where $m$ is odd, we construct~$G$ on $n$ vertices as follows.
Partition the vertices into $4$ parts $A,B,C,D$, with $|A|=|C|=m$, $|B|=m+1$ and $|D|=m+2$.
Each of $A$ and $C$ spans a regular tournament, $B$ and $D$ are joined by a bipartite tournament
(i.e.~an orientation of the complete bipartite graph) which is as regular as possible.
We also add all edges from $A$ to $B$, from $B$ to $C$, from $C$ to $D$ and from $D$ to $A$.
Since every path which joins two vertices in~$D$ has to pass through~$B$, it follows that
every cycle contains at least as many vertices from~$B$ as it contains from~$D$.
As $|D|>|B|$ this means that one cannot cover all the vertices of~$G$ by disjoint cycles.
%i.e.~$G$ does not contain a 1-factor.
This construction can be extended to arbitrary~$n$
(see~\cite{KKOexact}). The following result exactly matches this bound and improves earlier ones of
several authors, e.g.~\cite{HaggkvistHamilton,HaggkvistThomasonHamilton,kellyKO}. 
\begin{theorem}{\bf (Keevash, K\"uhn and Osthus~\cite{KKOexact})} \label{main}
There exists an integer $n_0$ so that any oriented graph
$G$  on $n \ge n_0$ vertices with minimum
semidegree $\delta^0(G) \ge \frac{3n-4}{8}$ contains a Hamilton cycle.
\end{theorem}
The proof of this result is based on some ideas in~\cite{kellyKO}.
H\"aggkvist~\cite{HaggkvistHamilton} also made the following conjecture which is closely related
to Theorem~\ref{main}. Given an oriented graph~$G$, let~$\delta(G)$ denote the minimum degree of~$G$
(i.e.~the minimum number of edges incident to a vertex) and set
$\delta^*(G):=\delta(G)+\delta^+(G)+\delta^-(G)$.
\begin{conj}{\bf (H\"aggkvist~\cite{HaggkvistHamilton})} \label{haggconj}
Every oriented graph~$G$ on $n$ vertices with $\delta^*(G)>(3n-3)/2$ contains a Hamilton cycle.
\end{conj}
(Note that this conjecture does not quite imply Theorem~\ref{main} as it results in a marginally
greater minimum semidegree condition.)
In~\cite{kellyKO}, Conjecture~\ref{haggconj} was verified approximately, i.e.~if
$\delta^*(G) \ge (3/2+o(1))n$, then $G$ has a Hamilton cycle (note this implies an
approximate version of Theorem~\ref{main}).
The same methods also yield an approximate version of Theorem~\ref{woodall} for oriented graphs.
\begin{theorem}{\bf (Kelly, K\"uhn and Osthus~\cite{kellyKO})}\label{thm:Ore}
For every $\alpha>0$ there exists an integer $n_0=n_0(\alpha)$ such that every oriented
graph~$G$ of order $n\geq n_0$ with $d^+(x)+d^-(y)\ge (3/4+\alpha)n$ whenever $G$ does
not contain an edge from~$x$ to~$y$ contains a Hamilton cycle.
\end{theorem}
The above construction of H\"aggkvist shows that the bound is best possible up to the term 
$\alpha n$. It would be interesting to obtain an exact version of this result.

Note that Theorem~\ref{main} implies that every sufficiently large regular
tournament on~$n$ vertices contains at least $n/8$ edge-disjoint Hamilton cycles.
(To verify this, note that in a regular tournament, all in- and outdegrees are equal to $(n-1)/2$.
We can then greedily remove Hamilton cycles
as long as the degrees satisfy the condition in Theorem~\ref{main}.)
It is the best bound so far towards the following conjecture of Kelly~(see e.g.~\cite{digraphsbook}).
\begin{conj}{\bf (Kelly)} 
Every regular tournament on~$n$ vertices can be
partitioned into~$(n-1)/2$ edge-disjoint Hamilton cycles.
\end{conj}
A result of Frieze and Krivelevich~\cite{FKhampack} states that 
every dense $\eps$-regular digraph contains a collection of edge-disjoint Hamilton cycles which covers
almost all of its edges. This implies that the same holds for almost every 
tournament. Together with a lower bound by McKay~\cite{McKay} on the number of 
regular tournaments, it is easy to see that 
the above result in~\cite{FKhampack} also implies that almost every regular tournament
contains a collection of edge-disjoint Hamilton cycles which covers
almost all of its edges.
Thomassen made the following conjecture which replaces the assumption of regularity by
high connectivity. 
\begin{conj}{\bf (Thomassen~\cite{thomassenconj})}
For every $k \ge 2$ there is an integer $f(k)$ so that every strongly $f(k)$-connected
tournament has $k$ edge-disjoint Hamilton cycles.
\end{conj}
The following conjecture of Jackson is also closely related to Theorem~\ref{main} -- it
would imply a much better degree condition for regular oriented graphs.
\begin{conj}{\bf (Jackson~\cite{jacksonconj})} \label{jacksonconj}
For $d>2$, every $d$-regular oriented graph $G$ on $n \le 4d+1$ vertices 
is Hamiltonian.
\end{conj}
The disjoint union of two regular tournaments on $n/2$ vertices shows that this would be best possible.
An undirected analogue of Conjecture~\ref{jacksonconj} was proved by Jackson~\cite{jacksonreg}.
It is easy to see that every tournament on~$n$ vertices with minimum semidegree at least
$n/4$ has a Hamilton cycle. In fact, for tournaments $T$ of large order~$n$ with minimum semidegree at least
$n/4+\eps n$, Bollob\'as and H\"aggkvist~\cite{BHpower} 
proved the stronger result that (for fixed $k$) $T$ even contains the 
$k$th power of a Hamilton cycle. It would be interesting to 
find corresponding degree conditions which ensure this for arbitrary digraphs and for
oriented graphs.

\subsection{Degree sequences forcing Hamilton cycles in directed graphs}

For undirected graphs, Dirac's theorem is generalized by Chv\'atal's
theorem~\cite{chvatal} that characterizes all those degree sequences which
ensure the existence of a Hamilton cycle in a graph: suppose that the degrees of the graph
are $d_1\le \dots \le d_n$. If $n \geq 3$ and $d_i \geq i+1$ or $d_{n-i} \geq n-i$
for all $i <n/2$ then $G$ is Hamiltonian. This condition on the degree sequence is
best possible in the sense that for any degree sequence violating this condition there
is a corresponding graph with no Hamilton cycle.
Nash-Williams~\cite{nw}
raised the question of a digraph analogue of Chv\'atal's theorem quite soon after the
latter was proved:
for a digraph~$G$ it is natural to consider both its outdegree sequence $d^+ _1,\dots , d^+ _n$
and its indegree sequence $d^- _1,\dots , d^- _n$. Throughout, we take the convention that
$d^+ _1\le \dots \le d^+ _n$ and $d^- _1 \le \dots \le  d^- _n$ without mentioning this explicitly.
Note that the terms $d^+ _i$ and $d^- _i$ do not necessarily correspond to the degree of the same vertex
of~$G$. 
\begin{conj}[Nash-Williams~\cite{nw}]\label{nw}
Suppose that $G$ is a strongly connected digraph on $n \geq 3$ vertices
such that for all $i < n/2$
\begin{itemize}
\item[{\rm (i)}]  $d^+ _i \geq i+1 $ or $ d^- _{n-i} \geq n-i $,
\item[{\rm (ii)}] $ d^- _i \geq i+1$ or $ d^+ _{n-i} \geq n-i.$
\end{itemize}
Then $G$ contains a Hamilton cycle.
\end{conj}
It is even an open problem whether the conditions imply the existence of a cycle through
any pair of given vertices (see~\cite{bt}). 
It is easy to see that one cannot omit the condition that~$G$ is strongly connected. 
The following example (which is a straightforward generalization of the corresponding 
undirected example) shows that the 
degree condition in Conjecture~\ref{nw} would be best possible in the sense that for all
$n\ge 3$ and all $k<n/2$ there is a non-Hamiltonian strongly connected digraph~$G$ on~$n$ vertices
which satisfies the degree conditions except that $d^+_k,d^-_k\ge k+1$ are replaced by
$d^+_k,d^-_k\ge k$ in the $k$th pair of conditions.
To see this, take an independent set~$I$ of size $k<n/2$ 
and a complete digraph~$K$ of order~$n-k$. Pick a set~$X$ of~$k$ vertices of~$K$
and add all possible edges (in both directions) between~$I$ and~$X$. The digraph~$G$
thus obtained is strongly connected, not Hamiltonian and
$$\underbrace{k, \dots ,k}_{k \text{ times}}, \underbrace{n-1-k, \dots , n-1-k}_{n-2k
\text{ times}}, \underbrace{n-1, \dots , n-1}_{k \text{ times}}$$ is both the out- and
indegree sequence of~$G$. In contrast to the undirected case there exist examples with 
a similar degree sequence to the above but whose structure is quite different (see~\cite{KOTchvatal}).
In~\cite{KOTchvatal}, the following approximate version of Conjecture~\ref{nw} for large digraphs was proved. 

\begin{theorem}[K\"uhn, Osthus and Treglown \cite{KOTchvatal}]\label{approxnw}
For every $\alpha >0$ there exists an integer $n_0 =n_0 (\alpha)$ such that the following holds.
Suppose $G$ is a digraph on $n \geq n_0$ vertices  such that for all $i < n/2$
\begin{itemize}
\item $ d^+ _i \geq i+ \alpha n  $ or $ d^- _{n-i- \alpha n} \geq n-i $,
\item $ d^- _i \geq i+ \alpha n $ or $ d^+ _{n-i- \alpha n} \geq n-i .$
\end{itemize}Then $G$ contains a Hamilton cycle.
\end{theorem}
Theorem~\ref{approxnw} was derived from a result in~\cite{KKOexact} on the existence of a 
Hamilton cycle in an oriented graph satisfying a `robust' expansion property.

The following weakening of Conjecture~\ref{nw} was posed earlier by Nash-Williams \cite{ch2}.
It would yield a digraph analogue of P\'osa's theorem which states that a graph~$G$ on~$n\ge 3$
vertices has a Hamilton cycle if its degree sequence $d_1\le  \dots \le d_n$ satisfies $d_i \geq i+1$ 
for all $i<(n-1)/2$ and if additionally $d_{\lceil n/2\rceil} \geq \lceil n/2\rceil$ when~$n$ is odd~\cite{posa}. 
Note that P\'osa's theorem is much stronger than Dirac's theorem but is a special case of Chv\'atal's theorem.
\begin{conj}[Nash-Williams~\cite{ch2}]\label{nw2}
Let $G$ be a digraph on $n \geq 3$ vertices such that $d^+ _i,d^-_i \geq i+1 $
for all $i <(n-1)/2$ and such that additionally
$d^+_{\lceil n/2\rceil},d^-_{\lceil n/2\rceil} \geq \lceil n/2\rceil$ when~$n$ is odd. 
Then~$G$ contains a Hamilton cycle.
\end{conj}
The previous example shows the degree condition would be best possible in the same sense as described there. 
The assumption of strong connectivity is not necessary in Conjecture~\ref{nw2},
as it follows from the degree conditions.
Theorem~\ref{approxnw} immediately implies an approximate version of Conjecture~\ref{nw2}.

It turns out that the conditions of Theorem~\ref{approxnw} even
guarantee the digraph~$G$ to be \emph{pancyclic}, i.e.~$G$ contains a cycle of length~$t$
for all $t=2,\dots,n$.
Thomassen~\cite{tom} as well as H\"aggkvist and Thomassen~\cite{haggtom}
gave degree conditions which imply that every digraph with  minimum
semidegree $>n/2$ is pancyclic. The latter  bound can also be deduced directly from 
Theorem~\ref{hamGH}.
The complete bipartite digraph whose vertex class sizes
are as equal as possible shows that the bound is best possible.
For oriented graphs the minimum semidegree threshold which guarantees pancyclicity
turns out to be $(3n-4)/8$ (see~\cite{kellyKOpan}).

\subsection{Powers of Hamilton cycles in graphs}

The following result is a common extension (for large $n$) of Dirac's theorem and the 
Hajnal-Szemer\'edi theorem. It was originally conjectured (for all $n$) by Seymour. 
\begin{theorem} \label{KSSpowers} {\bf (Koml\'os, S\'ark\"ozy and Szemer\'edi~\cite{KSSz98})}
For every $k\ge 1$ there is an integer $n_0$ so that every graph $G$ on $n \ge n_0$ vertices and with
$\delta(G) \ge \frac{k}{k+1} n$ contains the $k$th power of a Hamilton cycle.
\end{theorem}
Complete $(k+1)$-partite graphs whose vertex classes have almost (but not exactly) equal size
show that the minimum degree bound is best possible. 
Previous to this a large number of partial results had been proved (see~e.g.~\cite{LSS} for a history of the problem).
Very recently, Levitt, Sark\"ozy and Szemer\'edi~\cite{LSS} gave a proof of the case $k=2$ 
which avoids the use of the Regularity lemma, resulting in a much better bound on~$n_0$.
Their proof is based on a technique introduced by R\"odl, Ruci\'nski and Szemer\'edi~\cite{RRSDirac}
for hypergraphs.
The idea of this method (as applied in~\cite{LSS}) is first to find an `absorbing' path $P^2$:
roughly, $P^2$ is the second power of a path $P$ which, given any vertex $x$, has the property
that $x$ can be inserted into $P$ so that $P \cup {x}$
still induces the second power of a path. The proof of the existence of $P^2$ is heavily based
on probabilistic arguments. Then one finds the second power $Q^2$ of a path which is almost spanning in $G-P^2$.
One can achieve this by repeated applications of the Erd\H{o}s-Stone theorem.
One then connects up~$Q^2$ and~$P^2$ into the second power of a cycle and finally uses the absorbing 
property of $P^2$ to incorporate the vertices left over so far.

\subsection{Hamilton cycles in hypergraphs} \label{hypercycle}

It is natural to ask whether one can generalize Dirac's theorem to uniform hypergraphs.
There are several possible notions of a hypergraph cycle. One generalization of the 
definition of a cycle in a graph is the following one.
An $r$-uniform hypergraph $C$ is a \emph{cycle of order $n$}
if there a exists a cyclic ordering $v_1,\dots,v_n$ of its~$n$ vertices such 
that every consecutive pair $v_iv_{i+1}$ lies in a hyperedge of $C$
and such that every hyperedge of $C$ consists of consecutive 
vertices. Thus the cyclic ordering of the vertices of $C$ induces a
cyclic ordering of its hyperedges. A cycle is \emph{tight} if every
$r$ consecutive vertices form a hyperedge. 
A cycle of order $n$ is \emph{loose} if all pairs of consecutive edges (except possibly one pair)
have exactly one vertex in common. (So every tight cycle contains a spanning
loose cycle but a cycle might not necessarily contain a spanning loose cycle.)
There is also the even more general notion 
of a \emph{Berge-cycle}, which consists of a sequence of vertices where each pair of consecutive 
vertices is contained in a common hyperedge.
\noindent
\begin{center}
\psfrag{tight cycle}[][]{\small{tight cycle}}
\includegraphics[scale=0.3]{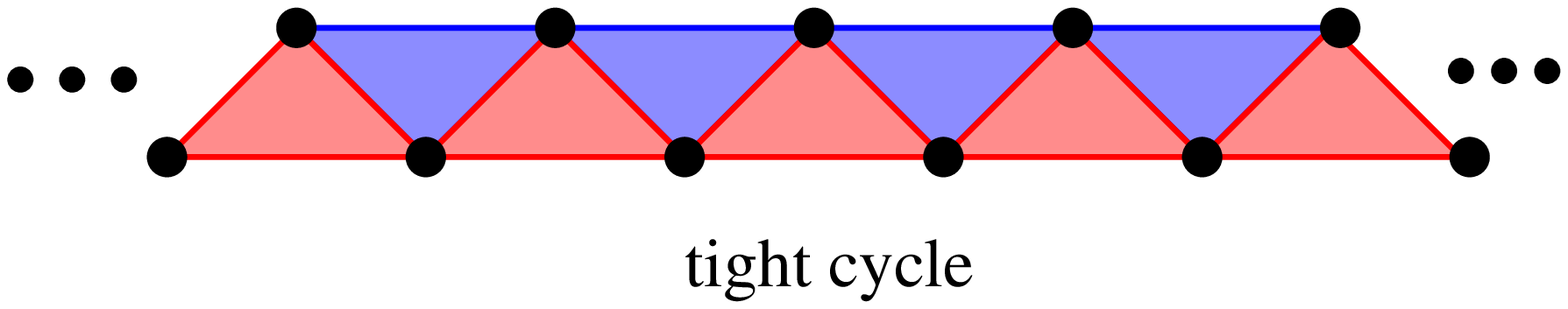}
\end{center}
\begin{center}
\psfrag{cycle}[][]{\small{cycle}}
\includegraphics[scale=0.3]{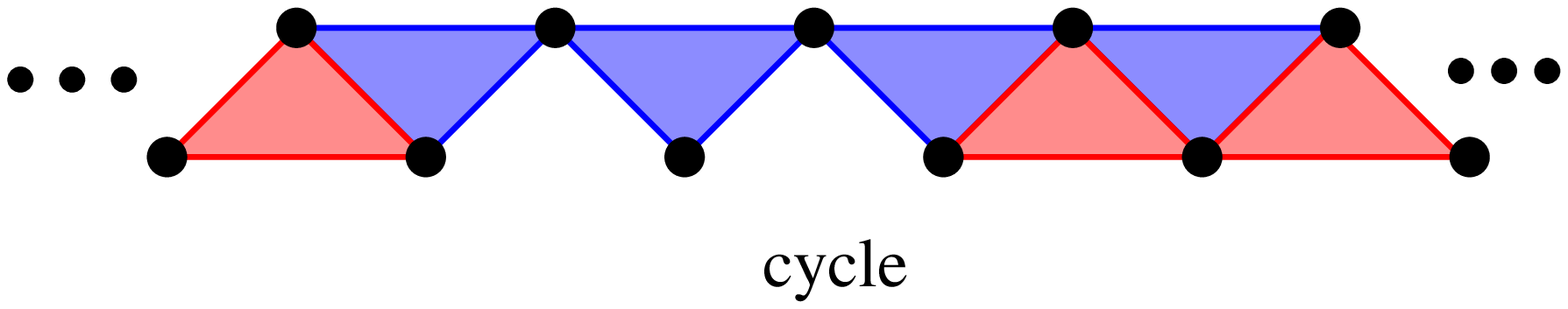}
\end{center}
\begin{center}
\psfrag{loose cycle of even order}[][]{\small{loose cycle}}
\includegraphics[scale=0.3]{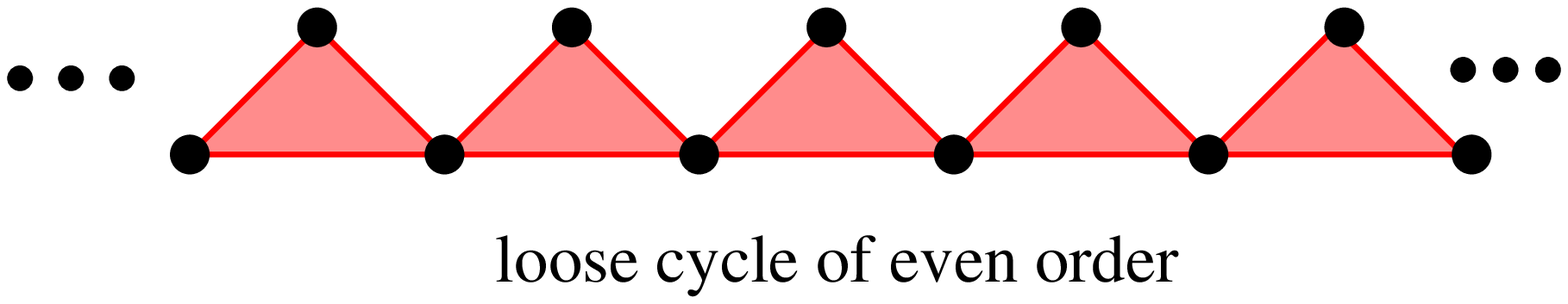}
\end{center}
A \emph{Hamilton cycle} of a uniform hypergraph $G$ is a subhypergraph of
$G$ which is a cycle containing all its vertices. Theorem~\ref{tightcycle} gives an analogue
of Dirac's theorem for tight hypergraph cycles, while Theorem~\ref{loosecycle} gives an analogue for
$3$-uniform (loose) cycles.
\begin{theorem}{\bf (R\"odl, Ruci\'nski and Szemer\'edi~\cite{RRSDirac})} \label{tightcycle}
For all $r \in \mathbb{N}$ and $\alpha>0$ there is an integer $n_0=n_0(r,\alpha)$ such that every
$r$-uniform hypergraph $G$ with $n\ge n_0$ vertices and minimum degree 
at least $n/2+\alpha n$ contains a tight Hamilton cycle.
\end{theorem}\begin{theorem}{\bf (Han and Schacht~\cite{HScycle};
Keevash, K\"uhn, Mycroft and Osthus~\cite{rcycle})} \label{loosecycle}
For all $r \in \mathbb{N}$ and  $\alpha>0$ there is an integer $n_0=n_0(\alpha)$ such that every
$r$-uniform hypergraph $G$ with $n\ge n_0$ vertices and minimum degree 
at least $n/(2r-2)+\alpha n$ contains a loose Hamilton cycle.
\end{theorem}
Both results are best possible up to the error term $\alpha n$. In fact,
if the minimum degree is less than $\lceil n/(2r-2) \rceil$, 
then we cannot even guarantee \emph{any} Hamilton cycle in an $r$-uniform hypergraph.
The case $r=3$ of Theorems~\ref{tightcycle} and~\ref{loosecycle} was proved earlier 
in~\cite{RRSDirac3} and~\cite{KOloose} respectively.
The result in~\cite{HScycle} also covers the notion of an $r$-uniform $\ell$-cycle
for $\ell<r/2$ (here we ask for consecutive edges to intersect in precisely $\ell$ vertices).
Hamiltonian Berge-cycles were considered by Bermond et al.~\cite{bermond}.
%%%%%%%%%%%%%%%%%%%%%%%%%%%%%%%%%%%%%%%%%%%%%%%%%%%%%%%%%%%%%%%%%%%%%%%%%%%%%%%%%%%%

\section{Bounded degree spanning subgraphs}

Bollob\'as and Eldridge~\cite{BE78} as well as Catlin~\cite{catlin} made the
following very general conjecture on embedding graphs.
If true, this conjecture would be a far-reaching generalization of the Hajnal-Szemer\'edi theorem
(Theorem~\ref{hajnalsz}).
\begin{conj}[Bollob\'as and Eldridge~\cite{BE78}, Catlin~\cite{catlin}]\label{conjBE}
If $G$ is a graph on $n$ vertices with $\delta(G) \ge \frac{\Delta n-1}{\Delta+1}$, then
$G$ contains any graph $H$ on $n$ vertices with maximum degree at most $\Delta$.
\end{conj}
The conjecture has been proved for graphs $H$ of maximum degree at most 2~\cite{AB93,AF96} and for
large graphs of maximum degree at most 3~\cite{CSS03}.
Recently, Csaba~\cite{Csababipartite} proved it for bipartite graphs $H$ of arbitrary 
maximum degree $\Delta$, provided the order of~$H$ is sufficiently large compared to $\Delta$. 
In many applications of the Blow-up lemma, the graph~$H$ is embedded into~$G$ by splitting~$H$
up into several suitable parts and applying the Blow-up lemma to each of these parts
(see e.g.~the example in Section~\ref{sample}). It is not clear how to achieve this
for~$H$ as in Conjecture~\ref{conjBE}, as $H$ may be an `expander'.
So the proofs in~\cite{Csababipartite,CSS03} rely on a variant of the Blow-up lemma which is suitable
for embedding such `expander graphs'. Also, Kaul, Kostochka and Yu~\cite{KKY} showed
(without using the Regularity lemma) that the conjecture holds if we increase the minimum degree condition 
to $\frac{\Delta n+2n/5-1}{\Delta+1}$.

Theorem~\ref{KSS} suggests that one might replace $\Delta$ in Conjecture~\ref{conjBE} with
$\chi(H)-1$, resulting in a smaller minimum degree bound for some graphs $H$.
This is far from being true in general (e.g.~let $H$ be a $3$-regular bipartite expander and let~$G$ be the union
of two cliques which have equal size and are almost disjoint).
However, Bollob\'as and Koml\'os conjectured that this does turn out to be true 
if we restrict our attention to a certain class of `non-expanding' graphs. This conjecture was recently confirmed
in~\cite{bandwidth}. The bipartite case was proved earlier by Abbasi~\cite{Abbasi}.
\begin{theorem}{\bf (B\"ottcher, Schacht and Taraz~\cite{bandwidth})} \label{bandwidth}
For every $\gamma>0$ and all integers $r \ge 2$ and $\Delta$, there exist $\beta>0$ and $n_0$ with the 
following property.
Every  graph $G$ of order $n \ge n_0$ and minimum degree at least
$(1-1/r+\gamma)n$ contains every $r$-chromatic graph $H$ of order $n$, maximum degree at most $\Delta$ and
bandwidth at most $\beta n$ as a subgraph.
\end{theorem}
Here the \emph{bandwidth} of a graph $H$ is the smallest integer
$b$ for which there exists an enumeration $v_1,\dots,v_{|H|}$ of the vertices
of $H$ such that every edge $v_iv_j$ of $H$ satisfies $|i-j|\le b$.
Note that $k$th powers of cycles have bandwidth
$2k$, so Theorem~\ref{bandwidth} implies an approximate version of Theorem~\ref{KSSpowers}. 
(Actually, this is only the case if $n$ is a multiple of $k+1$, as otherwise
the $k$th power of a Hamilton cycle fails to be $(k+1)$-colourable. 
But~\cite{bandwidth} contains a more general result which allows for a small number of
vertices of colour $k+2$.)
A further class of graphs having small bandwidth and bounded degree are 
planar graphs with bounded degree~\cite{bandwidthplanar}.
(See~\cite{KOtriang,KOTplanar} for further results on embedding planar graphs in graphs
of large minimum degree.)
Note that the discussion in Section~\ref{packings} implies that
the minimum degree bound in Theorem~\ref{bandwidth} 
is approximately best possible for certain graphs $H$ but not for all graphs.
Abbasi~\cite{abbasiband} showed that there are graphs $H$ for which the linear error term $\gamma n$ in 
Theorem~\ref{bandwidth} is necessary. One might think that one could
reduce the error term  to a constant for graphs of bounded bandwidth. 
However, this turns out to be incorrect. (We grateful to Peter Allen for pointing this out to us.) 

%(So this would apply
%to e.g.~bounded powers of cycles but not to the square grid on $n$ vertices, 
%which has bandwidth $\Theta(\sqrt{n})$.) 
% 
%\begin{conj} \label{bandwidthconst}
%For all integers $b$ and $r$ there exists a
%constant $C=C(b,r)$ such that each graph of order $n$ and minimum degree at least
%$(1-1/r)n+C$ contains every $r$-chromatic graph $H$ of order $n$ and
%bandwidth at most $b$ as a subgraph.
%\end{conj}
%Note that Conjecture~\ref{bandwidthconst} would be a generalization of Theorem~\ref{KSS}.
Alternatively, one can try to replace the bandwidth assumption in Theorem~\ref{bandwidth}
with a less restrictive parameter. For instance, Csaba~\cite{csabasep} gave a minimum degree
condition on $G$ which guarantees a copy of a `well-separated' graph $H$ in $G$.
Here a graph with $n$ vertices is \emph{$\alpha$-separable} if there is a set $S$ of vertices of 
size at most $\alpha n$ so that all components of $H-S$ have size at most $\alpha n$.
It is easy to see that every graph with $n$ vertices and bandwidth at most $\beta n$
is $\sqrt{\beta}$-separable. (Moreover large trees are $\alpha$-separable for 
$\alpha \to 0$ but need not have small bandwidth, so considering separability is less restrictive than
bandwidth.)
%A further way of generalizing Theorem~\ref{bandwidth} would be to be replace bandwidth by a less
%restrictive parameter. For instance, Csaba introduced the notion of a `well-separated' graph:
%a graph with $n$ vertices is \emph{$\alpha$-separable} if there is a set $S$ of vertices of 
%size at most $\alpha n$ so that all components of $H-S$ have size at most $\alpha n$.
%It is easy to see that every graph with $n$ vertices and bandwidth at most $\beta n$
%is $\sqrt{\beta}$-separable. (Moreover large trees are $\alpha$-separable for 
%$\alpha \to 0$ but need not have small bandwidth, so considering separability is less restrictive than
%bandwidth.)
%In particular, the following result implies the bipartite case of Theorem~\ref{bandwidth}.
%\begin{theorem}{\bf (Csaba~\cite{csabaseparate})} \label{csabasep}
%For every $\gamma>0$ and all integers $\Delta$ and $r\ge 2$ there exist $\beta>0$ and $n_0$
%with the following property: every graph $G$ of order $n \ge n_0$ and minimum degree at least
%$(1-1/(2(r-1))+\gamma)n$ contains every $r$-chromatic $\beta$-separable graph $H$ of order $n$
%with maximum degree at most~$\Delta$ as a subgraph.
%\end{theorem}
%It would be interesting to know whether one can replace the term $1/2(r-1)$ by $1/r$ and thus obtain a 
%generalization of Theorem~\ref{bandwidth}.

Here is another common generalization of Dirac's theorem and the triangle case of Theorem~\ref{hajnalsz}
(i.e.~the Corr\'adi-Hajnal theorem).
It proves a 
conjecture by El-Zahar (actually, El-Zahar made the conjecture for all values of $n$, this is
still open).%
\COMMENT{for future proposals: can we use absorbing technique to get new proof?}
\begin{theorem} {\bf (Abbasi~\cite{Abbasi})} \label{abbasi}
There exists an integer $n_0$ so that the following holds. Suppose that $G$ is a graph on $n \ge n_0$ vertices 
and $n_1,\dots, n_k\ge 3$ are so that 
$$\sum_{i=1}^k n_i=n \qquad and \qquad \delta(G) \ge \sum_{i=1}^k \lceil n_i/2 \rceil.
$$
Then $G$ has $k$ vertex-disjoint cycles whose lengths are $n_1,\dots,n_k$.
\end{theorem}
%Given $n_i$ as above, let $\ell$ denote the number of cycle lengths which are odd. 
%A complete $3$-partite graph on $n$ vertices with one vertex class of size $\ell-1$ and where
%the size of the two others is as equal as possible shows that the bound is best possible. 
%Note that if $H$ is a union of $k$ vertex disjoint cycles as above, then 
%$\chi_{cr}(H)=2n/(n-\# odd)$, where $\# odd$ denotes the number of cycles whose length is odd.
%It follows that $\sum_{i=1}^k \lceil n_i/2 \rceil=(1-1/\chi_{cr}(H))n$.
Note that $\sum_{i=1}^k \lceil n_i/2 \rceil=\sum_{i=1}^k(1-1/\chi_{cr}(C_i))n_i$,
where $C_i$ denotes a cycle of length $n_i$.
This suggests the following more general question (which was raised by Koml\'os \cite{JKtiling}):
Given $t\in\mathbb{N}$, does there exists an $n_0=n_0(t)$ such that whenever
$H_1,\dots,H_k$ are graphs which each have at most $t$ vertices and which
together have $n\ge n_0$ vertices and whenever~$G$ is a graph on $n$
vertices with minimum degree at least $\sum_i (1-1/\chi_{cr}(H_i))|H_i|$,
then there is a set of vertex-disjoint copies of $H_1,\dots,H_k$ in $G$?
In this form, the question has a negative answer by (the lower bound in) 
Theorem~\ref{thmmaingeneral}, but it would be interesting to find a common generalization of
Theorems~\ref{thmmaingeneral} and~\ref{abbasi}.

It is also natural to ask corresponding questions for oriented and directed graphs.
As in the case of Hamilton cycles, the questions appear much harder than in the undirected case 
and again much less is known. 
Keevash and Sudakov~\cite{keevashsudakov} recently obtained the following result which can be viewed
as an oriented version of the $\Delta=2$ case of Conjecture~\ref{conjBE}.
\begin{theorem} {\bf (Keevash and Sudakov~\cite{keevashsudakov})}
There exist constants $c,C$ and an integer~$n_0$ so that whenever $G$ is an oriented graph on $n\ge n_0$ vertices
with minimum semidegree at least $(1/2-c)n$ and whenever $n_1,\dots,n_t$ are so that 
$\sum_{i=1}^t n_i \le n-C$, then $G$ contains disjoint cycles of length $n_1,\dots,n_t$.
\end{theorem}
In the case of triangles (i.e.~when all the $n_i=3$), they show that one can choose $C=3$
(one cannot take $C=0$).
\cite{keevashsudakov} also contains a discussion of related open questions for tournaments
and directed graphs.
Similar questions were also raised earlier by Song~\cite{song}. For instance,
given $t$, what is the smallest integer $f(t)$ so that all but a finite number of $f(t)$-connected
tournaments $T$ satisfy the following: Let $n$ be the number of vertices of $T$ and let 
$\sum_{i=1}^t n_i = n$. Then $T$ contains disjoint cycles of length $n_1,\dots,n_t$.
%%%%%%%%%%%%%%%%%%%%%%%%%%%%%%%%%%%%%%%%%%%%%%%%%%%%%%%%%%
\section{Ramsey Theory}
The Regularity lemma can often be used to show that the Ramsey numbers of sparse graphs $H$ are small.
(The \emph{Ramsey number $R(H)$} of~$H$ is the smallest $N \in \mathbb{N}$ such that for every
$2$-colouring of the complete graph on $N$ vertices one can find a monochromatic copy of~$H$.)
In fact, the first result which demonstrated the use of the Regularity lemma in extremal graph
theory was the following result of Chv\'atal, R\"odl, Szemer\'edi and Trotter~\cite{CRST},
which states that graphs of bounded degree have linear Ramsey numbers: 
\begin{theorem}{\bf (Chv\'atal, R\"odl, Szemer\'edi and Trotter~\cite{CRST})} \label{CRST}
For all $\Delta\in\mathbb{N}$ there is a constant $C=C(\Delta)$ so that every graph $H$ with maximum degree 
$\Delta(H) \le \Delta$ and $n$ vertices satisfies $R(H)\le Cn$.
\end{theorem}
The constant $C$ arising from the original proof (based on the Regularity lemma) is quite large.
The bound was improved in a series of papers. Recently, Fox and Sudakov~\cite{FoxSudakov} showed that 
$R(H)\le 2^{4\chi(H)\Delta}n$ (the bipartite case was also proved 
independently by Conlon~\cite{ConlonRamsey}). For bipartite graphs, a construction
from~\cite{GRR} shows that this bound is best possible apart from the value of the absolute constant~$4\cdot 2$ 
appearing in the exponent. %

Theorem~\ref{CRST} was recently generalized to hypergraphs~\cite{CFKO3,CFKOk,NORS,Ishigami} using hypergraph
versions of the Regularity lemma. Subsequently, Conlon, Fox and Sudakov~\cite{CFSRamsey} gave a shorter proof 
which gives a better constant and does not rely on the Regularity lemma. 

One of the most famous conjectures in Ramsey theory is the Burr-Erd\H{o}s conjecture on $d$-degenerate graphs, 
which generalizes Theorem~\ref{CRST}. Here a graph~$G$ is \emph{$d$-degenerate} if every subgraph has a vertex of
degree at most~$d$. In other words, $G$ has no `dense' subgraphs.
\begin{conj}{\bf (Burr and Erd\H{o}s~\cite{BERamsey})}\label{burrerdos}
For every $d$ there is a constant $C=C(d)$ so that every $d$-degenerate graph $H$ on $n$ vertices
satisfies $R(H)\le Cn$.
\end{conj}
It has been proved in many special cases (see e.g.~the introduction of~\cite{FoxSudakovbip} 
for a recent overview).
Also, Kostochka and Sudakov~\cite{KSRamsey} proved that it is `approximately' true:
\begin{theorem}{\bf (Kostochka and Sudakov~\cite{KSRamsey})} \label{approxburrerdos}
For every $d$ there is a constant $C=C(d)$ so that every $d$-degenerate graph $H$ on $n$ vertices 
satisfies $R(H)\le 2^{C (\log n)^{2d/(2d+1)}} n $.
\end{theorem}
The exponent `$2d/(2d+1)$' of the logarithm was improved to `1/2' in~\cite{FoxSudakovbip}. All the results
in~\cite{ConlonRamsey,FoxSudakov,FoxSudakovbip,KSRamsey} rely on variants of the same 
probabilistic argument, which was first applied to special cases of Conjecture~\ref{burrerdos} in~\cite{KoRodl}. 
To give an idea of this beautiful argument, we use a simple version to give a proof of the
following density result (which is implicit in several of the above papers):
it implies that bipartite graphs $H$ whose maximum degree is logarithmic in their order
have polynomial Ramsey numbers. (The logarithms in the statement and the proof are binary.)
\begin{theorem} \label{random}
Suppose that $H=(A',B',E')$ is a bipartite graph on $n \ge 2$ vertices and $\Delta(H) \le \log n$.
Suppose that $m \ge n^8$.
Then every bipartite graph $G=(A,B,E)$ with $|A|=|B|=m$ and at least $m^2/8$ edges contains a copy of $H$.
In particular, $R(H) \le 2n^8$.
\end{theorem}
An immediate corollary is that the Ramsey number of a $d$-dimensional cube $Q_d$ is
polynomial in its number $n=2^d$ of vertices (this fact was first observed in~\cite{Shi} based on
an argument similar to that in~\cite{KoRodl}). The best current bound of
$R(Q_d) \le d 2^{2d+5}$ is given in~\cite{FoxSudakov}.
Burr and Erd\H{o}s~\cite{BERamsey} conjectured that the bound should actually be linear in~$n=2^d$.
\proof
Write $\Delta:=\log n$.
Let $b_1,\dots,b_s$ be a sequence of $s:=2\Delta$ not necessarily distinct vertices of $B$,
chosen uniformly and independently at random and write $S:=\{b_1,\dots,b_s \}$.
Let $N(S)$ denote the set of common neighbours of vertices in~$S$.
Clearly, $S \subseteq N(a)$ for every $a \in N(S)$.
So Jensen's inequality implies that
\begin{align*}
\ex ( |N(S)| ) & = \sum_{a \in A} \pr ( a \in N(S) )
= \sum_{a \in A} \left( \frac{|N(a)|}{m} \right)^s
= \frac{\sum_{a \in A} ( d(a) )^s}{m^s} \\
& \ge \frac{m \left( \frac{\sum_{a \in A} d(a) }{m} \right)^s  }{m^s} 
\ge  \frac{m  \left( (m^2/8)/m \right)^s  }{m^s} = \frac{m}{8^s} 
\ge \frac{n^8}{n^6} =n^2.
\end{align*}
We say that a subset $W \subseteq A$ is \emph{bad} if it has size $\Delta$ and
its common neighbourhood $N(W)$ satisfies $|N(W)| <n$.
Now let $Z$ denote the number of bad subsets $W$ of $N(S)$.
Note that the probability that a given set $W \subseteq A$ lies in $N(S)$ equals
$(|N(W)|/m)^s$ (since the probability that it lies in the neighbourhood of a fixed 
vertex $b\in B$ is $|N(W)|/m$). So 
$$
\ex Z= \sum_{W \text{bad}} \pr (W \subseteq N(S) )
\le \binom{m}{\Delta} \left( \frac{ n}{m} \right)^s 
\le m^\Delta  \left( \frac{ n}{m} \right)^s 
= \left( \frac{n^2}{m} \right)^\Delta \le (1/2)^\Delta <1.
$$
So $\ex ( |N(S)| -Z ) \ge n^2-1 \ge n$ and hence there is a choice of $S$ with $|N(S)| -Z \ge n$.
By definition, we can delete 
a vertex from every bad $W$ contained in $N(S)$ to obtain a set $T \subseteq N(S)$ 
with $|T| \ge n$
so that every subset $W \subseteq T$ with $|W|=\Delta$ satisfies $|N(W)| \ge n$.
Clearly we can now embed $H$: first embed $A'$ arbitrarily into $T$
and then embed the vertices of $B'$ one by one into $B$, 
using the property that $T$ has no bad subset.  

The bound on $R(H)$ can be derived as follows:
consider any $2$-colouring of the complete 
graph on $2n^8$ vertices. Partition its vertices arbitrarily into two sets $A$ and $B$ of size $n^8$
and then apply the main statement to the subgraph of~$G$ induced by the colour class 
having the most edges between $A$ and $B$.
\endproof
Note that the proof immediately shows that the bound on the maximum degree of $H$
can be relaxed: all we need is the property that 
every subgraph of $H$ has a vertex $b \in B'$ of low degree. 
In the proof of (the bipartite case) of
Theorem~\ref{approxburrerdos}, this was exploited as follows: roughly speaking one carries out the above 
argument twice (of course with different parameters than the above).
The first time we consider a random subset $S \subseteq B$ and the second time we consider a smaller random subset
$S' \subseteq T$. 

For some types of sparse graphs $H$, one can give even more precise 
estimates for $R(H)$ than the ones which follow from the above results. 
For instance, Theorem~\ref{zhao} has an immediate application to the Ramsey number of trees.
\begin{cor} \label{ramseytrees}
There is an integer $n_0$ so that if $T_n$ is a tree on $n \ge n_0$ vertices then $R(T_n) \le 2n-2$.
\end{cor}
Indeed, to derive Corollary~\ref{ramseytrees} from Theorem~\ref{zhao}, consider a $2$-colouring
of a complete graph $K_{2n-2}$ on $2n-2$ vertices, yielding a red graph $G_r$ and a blue graph $G_b$.
Order the vertices $x_i$ according to their degree (in ascending order) in $G_r$.
If $x_{n-1}$ has degree at least $n-1$ in $G_r$, then we can apply Theorem~\ref{zhao} to find
a red copy of $T$ in $G_r$.
If not, we can apply it to find a blue copy of $T$ in $G_b$.
For even $n$, the bound is best possible (let $T$ be a star and let $G_b$ and $G_r$ be regular of the same 
degree) and proves a conjecture of Burr and Erd\H{o}s~\cite{BEtrees}. For odd $n$, they conjectured
that the answer is $2n-3$.
Similarly, the Koml\'os-S\'os conjecture (Conjecture~\ref{KomlosSos}) would imply
that $R(T_n,T_m) \le n+m-2$, where $T_n$ and $T_m$ are trees on $n$ and $m$ vertices respectively.
Of course, Corollary~\ref{ramseytrees} is not best possible for every tree.
For instance, in the case when the tree is a path, Gerencs\'er and Gyarfas~\cite{GGpath} showed that
$R(P_n,P_n)= \lfloor (3n-2)/2 \rfloor$.
Further recent results on Ramsey numbers of paths and cycles (many of which 
rely on the Regularity lemma) can be found e.g.~in~\cite{3colourpath,FG}.
Hypergraph versions (i.e.~Ramsey numbers of tight cycles, loose cycles and Berge-cycles) were considered 
e.g.~in~\cite{ramseyhypercycle1,ramseyhypercycle2, bergeramsey}.%
    \COMMENT{
Luczak and FigaiWe show that for any real positive numbers $\alpha_1$, $\alpha_2$ and $\alpha_3$, the Ramsey number 
for a triple of even cycles of lengths $2\lfloor \alpha_1n\rfloor$, $2\lfloor \alpha_2n\rfloor$ and 
$2\lfloor \alpha_3n\rfloor$ is (asymptotically) equal to $(\alpha_1+\alpha_2+\alpha_3+\max\{\alpha_1,\alpha_2,\alpha_3\}+o(1))n$.
Also, Skokan has announced results on his homepage}

%%%%%%%%%%%%%%%%%%%%%%%%%%%%%%%%%%%%%%%%%%%%%%%%%%%%%%%%%%%%%%%%%%%%%%%%%%%%%%%%%%%%%%
\section{A sample application of the Regularity and Blow-up lemma} \label{sample}

In order to illustrate the details of the Regularity method for those not familiar with it, 
we now prove Theorem~\ref{KSS} for the
case when $H:=C_4$ and when we replace the constant~$C$ in the minimum degree condition
with a linear error term.

\begin{theorem}\label{thm:C4}
For every $0<\eta<1/2$ there exists an integer~$n_0$ such that every graph~$G$ whose order $n\ge n_0$
is divisible by~$4$ and whose minimum degree is at least~$n/2+\eta n$ contains a perfect
$C_4$-packing.
\end{theorem}

(Note that Theorem~\ref{thm:C4} also follows from Theorems~\ref{bandwidth} and~\ref{abbasi}.)
We start with the formal definition of $\eps$-regularity. 
The \emph{density} of a bipartite graph $G=(A,B)$ with vertex classes~$A$ and~$B$
is $$d_G(A,B):=\frac{e_G(A,B)}{|A||B|}.$$ We also write $d(A,B)$ if this is unambiguous.
Given $\eps>0$, we say that~$G$ is \emph{$\eps$-regular} if for all
sets $X\subseteq A$ and $Y\subseteq B$ with $|X|\ge \eps |A|$ and
$|Y|\ge \eps |B|$ we have $|d(A,B)-d(X,Y)|<\eps$. Given $d\in[0,1)$, we say
that $G$ is \emph{$(\eps,d)$-superregular} if all
sets $X\subseteq A$ and $Y\subseteq B$ with $|X|\ge \eps |A|$ and
$|Y|\ge \eps |B|$ satisfy $d(X,Y)>d$ and, furthermore, if $d_G(a)>d|B|$
for all $a\in A$ and $d_G(b)> d|A|$ for all $b\in B$.
Moreover, we will denote the neighbourhood of a vertex $x$ in a graph $G$ by $N_G(x)$.
Given disjoint sets~$A$ and~$B$ of vertices of~$G$, we write $(A,B)_G$ for the bipartite
subgraph of~$G$ whose vertex classes are~$A$ and~$B$ and whose edges are all the edges of~$G$
between~$A$ and~$B$.

Szemer\'edi's Regularity lemma~\cite{reglem} states that one can partition the vertices of
every large graph into a bounded number `clusters' so that most of the 
pairs of clusters induce $\eps$-regular bipartite graphs.
Proofs are also included in~\cite{BGraphTh} and~\cite{Diestel}.
Algorithmic proofs of the Regularity lemma were given in~\cite{algRL,FK}.
There are also several versions for hypergraphs (in fact, all the results in 
Section~\ref{hypercycle} are based on some hypergraph version of the Regularity lemma).
The first so-called `strong'  versions for $r$-uniform hypergraphs were proved in~\cite{gowers}
and~\cite{Count,RSkok}.  
\begin{lemma}[Szemer\'edi~\cite{reglem}] \label{reglem}
For all $\eps>0$ and all integers $k_0$ there is an $N=N(\eps,k_0)$ such
that for every graph~$G$ on~$n\ge N$ vertices there exists a
partition of $V(G)$ into $V_0,V_1,\dots,V_k$ such that the following holds:
\begin{enumerate}
\item[$\bullet$] $k_0\le k\le N$ and $|V_0|\le \eps n$,
\item[$\bullet$] $|V_1|=\dots=|V_k|=:m$,
\item[$\bullet$] for all but $\eps k^2$ pairs $1\le i<j\le k$ the graph $(V_i,V_j)_{G}$ is $\eps$-regular.
\end{enumerate}
\end{lemma}
Unfortunately, the constant $N$ appearing in the lemma is very large, 
Gowers~\cite{gowerstower} showed that it has at least a tower-type dependency on $\eps$.
We will use the
following degree form of Szemer\'edi's Regularity lemma which can be easily
derived from Lemma~\ref{reglem}. 
\begin{lemma}[Degree form of the Regularity lemma]\label{deg-reglemma}
For all $\eps>0$ and all integers $k_0$ there is an $N=N(\eps,k_0)$ such
that for every number $d\in [0,1)$ and for every graph~$G$ on~$n\ge N$ vertices there exist a
partition of $V(G)$ into $V_0,V_1,\dots,V_k$ and a spanning subgraph
$G'$ of $G$ such that the following holds:
\begin{enumerate}
\item[$\bullet$] $k_0\le k\le N$ and $|V_0|\le \eps n$,
\item[$\bullet$] $|V_1|=\dots=|V_k|=:m$,
\item[$\bullet$] $d_{G'}(x)>d_G(x)-(d+\eps)n$ for all vertices $x\in G$,
\item[$\bullet$] for all $i\ge 1$ the graph $G'[V_i]$ is empty,
\item[$\bullet$] for all $1\le i<j\le k$ the graph $(V_i,V_j)_{G'}$ is $\eps$-regular
and has density either $0$ or $>d$.
\end{enumerate}
\end{lemma}
The sets $V_i$ ($i\ge 1$) are called \emph{clusters}, $V_0$ is called
the \emph{exceptional set} and~$G'$ is called the \emph{pure graph}.

\medskip

\noindent
{\bf Sketch of proof of Lemma~\ref{deg-reglemma}\ }  
To obtain a partition as in Lemma~\ref{deg-reglemma}, apply Lemma~\ref{reglem}
with parameters $d,\eps', k_0'$ satisfying $1/k_0',\eps' \ll \eps,d,1/k_0$ to obtain
clusters $V'_1,\dots,V'_{k'}$ and an exceptional set~$V'_0$.
(Here $a \ll b < 1$ means that there is an increasing function $f$ such that
all the calculations in the argument work as long as $a \le f(b)$.)
Let $m':=|V'_1|=\dots=|V'_{k'}|$.
Now delete all edges between pairs of clusters which are not $\eps'$-regular and move any vertices
into $V'_0$ which were incident to at least $\eps n/10$ (say) of these deleted edges.
Secondly, delete all (remaining) edges between pairs of clusters whose density is at most $d+\eps'$.
Consider such a pair $(V'_i,V'_j)$ of clusters. For every vertex $x\in V'_i$ which has more than
$(d+2\eps')m'$ neighbours in~$V'_j$ mark all but $(d+2\eps')m'$ edges between~$x$ and~$V'_j$.
Do the same for the vertices in~$V'_j$ and more generally for all pairs of clusters of density at most $d+\eps'$.
It is easy to check that in total this yields at most $\eps'n^2$ marked edges. Move
all vertices into $V'_0$ which are incident to at least $\eps n/10$ of the marked edges.
Thirdly, delete any edges within the clusters.
Finally, we need to make sure that the clusters have equal size again (as we may have lost this
property during the deletion process). This can be done by splitting up the clusters into smaller
subclusters (which contain almost all the vertices and have equal size) and moving a small number of
further vertices into $V'_0$. A straightforward calculation shows that the new exceptional set
$V_0$ has size at most $\eps n$ as required.
\endproof

The \emph{reduced graph~$R$} is the graph whose vertices are
$1,\dots,k$ and in which~$i$ is joined to~$j$ whenever the bipartite subgraph
$(V_i,V_j)_{G'}$ of~$G'$ induced by~$V_i$ and~$V_j$ is $\eps$-regular and has density~$>d$.
Thus~$ij$ is an edge of~$R$ if and only if~$G'$ has an edge between~$V_i$ and~$V_j$.
Roughly speaking, the following result
states that~$R$ almost `inherits' the minimum degree of~$G$.

\begin{prop}\label{prop:mindegR}
If $0<2\eps\le d\le c/2$ and $\delta(G)\ge cn$ then $\delta(R)\ge (c-2d)|R|$.
\end{prop}
\proof
Consider any vertex~$i$ of~$R$ and pick $x\in V_i$. Then every neighbour of~$x$ in~$G'$
lies in $V_0\cup \bigcup_{j\in N_R(i)} V_j$. Thus $(c-(d+\eps))n\le d_{G'}(x)\le d_R(i)m +\eps n$
and so $d_R(i)\ge (c-2d)n/m\ge (c-2d)|R|$ as required.
\endproof

The proof of Proposition~\ref{prop:mindegR} is a point where it is important that~$R$ was defined
using the graph~$G'$ obtained from Lemma~\ref{deg-reglemma} and not using the partition
given by Lemma~\ref{reglem}.

In our proof of Theorem~\ref{thm:C4} the reduced graph~$R$ will contain a Hamilton path~$P$.
Recall that every edge~$ij$ of~$P\subseteq R$ corresponds to the $\eps$-regular bipartite
subgraph $(V_i,V_j)_{G'}$ of~$G'$ having density~$>d$. The next result shows that
by removing a small number of vertices from
each cluster (which will be added to the exceptional set~$V_0$) we can guarantee that the
edges of~$P$ even correspond to superregular pairs.

\begin{prop}\label{prop:superreg}
Suppose that $4\eps <d\le 1$ and that~$P$ is a Hamilton path in~$R$. Then every
cluster~$V_i$ contains a subcluster $V'_i\subseteq V_i$ of size~$m-2\eps m$ such
that $(V'_i,V'_j)_{G'}$ is $(2\eps,d-3\eps)$-superregular for every edge $ij\in P$.
\end{prop}
\proof
We may assume that $P=1\dots k$. Given any~$i< k$, the definition of regularity implies
that there are at most $\eps m$ vertices $x\in V_i$ such that $|N_{G'}(x)\cap V_{i+1}|\le (d-\eps)m$.
Similarly, for each $i>1$ there are at most $\eps m$ vertices $x\in V_i$ such that
$|N_{G'}(x)\cap V_{i-1}|\le (d-\eps)m$.
Let~$V'_i$ be a subset of size $m-2\eps m$ of~$V_i$ which contains none of the above vertices
(for all $i=1,\dots,k$). Then $V'_1,\dots, V'_k$ are as required.
\endproof

Of course, in Proposition~\ref{prop:superreg} it is not important that~$P$ is a Hamilton
path. One can prove an analogue whenever~$P$ is a subgraph of~$R$ of bounded maximum degree.
 We will also use the following special case of the Blow-up lemma of Koml\'os, S\'ark\"ozy and
Szemer\'edi~\cite{KSSblowup}. It implies that dense superregular pairs behave
like complete bipartite graphs with respect to containing bounded degree
graphs as subgraphs, i.e.~if the superregular pair has vertex classes~$V_i$ and~$V_j$
then any bounded degree bipartite graph on these vertex classes is a subgraph of this
superregular pair. An algorithmic version of the Blow-up lemma was proved by the same
authors in~\cite{algblowup}. A hypergraph version was recently proved by Keevash~\cite{keevash}.

\begin{lemma}[Blow-up lemma, bipartite case]\label{blowup}
Given $d>0$ and $\Delta\in \mathbb{N}$, there is a positive constant $\eps_0=\eps_0(d,\Delta)$
such that the following holds for every $\eps<\eps_0$. Given $m\in \mathbb{N}$,
let~$G^*$ be an $(\eps,d)$-superregular bipartite graph with vertex classes of size~$m$.
Then~$G^*$ contains a copy of every subgraph~$H$ of~$K_{m,m}$
with $\Delta(H)\le \Delta$.
\end{lemma}

\medskip

\noindent
{\bf Proof of Theorem~\ref{thm:C4}\ } 
We choose further positive constants~$\eps$ and~$d$
as well as $n_0\in\mathbb{N}$ such that
$$1/n_0\ll \eps\ll d\ll\eta<1/2.$$
(In order to simplify the exposition we will not determine these
constants explicitly.) We start by applying the degree form of the Regularity lemma (Lemma~\ref{deg-reglemma})
with parameters~$\eps$, $d$ and $k_0:=1/\eps$ to~$G$ to
obtain clusters $V_1,\dots,V_k$, an exceptional set~$V_0$, a pure graph~$G'$ and a reduced graph~$R$. 
Thus $k:=|R|$ and
\begin{equation}\label{eq:minR}
\delta(R)\ge (1/2+\eta -2d)k\ge  (1+\eta)k/2
\end{equation}
by Proposition~\ref{prop:mindegR}.
So~$R$ contains a Hamilton path~$P$ (this follows e.g.~from Dirac's theorem). By relabelling if
necessary we may assume that $P=1\dots k$. Apply Proposition~\ref{prop:superreg}
to obtain subclusters~$V'_i\subseteq V_i$ of size $m-2\eps m=:m'$ such that for every edge
$i(i+1)\in P$ the bipartite subgraph $(V'_i,V'_{i+1})_{G'}$ of~$G'$ induced by~$V'_i$ and~$V'_{i+1}$ is
$(2\eps,d/2)$-superegular. Note that the definition of $\eps$-regularity implies that
$(V'_i,V'_{j})_{G'}$ is still $2\eps$-regular of density
at least $d-\eps\ge d/2$ whenever $ij$ is an edge of~$R$.
We add all those vertices of~$G$ that are not contained in some~$V'_i$
to the exceptional set~$V_0$. Moreover, if~$k$ is odd then we also add all the vertices in~$V'_k$
to~$V_0$. We still denote the reduced graph by~$R$, its number of vertices by~$k$ and the
exceptional set by~$V_0$.
Thus%
    \COMMENT{Here we need that $m\le \eps n$, ie $n/k_0\le \eps n$. Ok if $k_0\ge 1/\eps$.}
$$ |V_0|\le \eps n+2\eps n+m \le 4\eps n.
$$
Let~$M$ denote the perfect matching in~$P$. So~$M$ consists of the edges $12,34,\dots, (k-1)k$.
The Blow-up lemma would imply that for every odd~$i$ the bipartite graph $(V'_i,V'_{i+1})_{G'}$
contains a perfect $C_4$-packing, provided that $2$ divides~$m'$. So we have already proved that~$G$
contains a $C_4$-packing covering \emph{almost} all of its vertices (this can also be easily
proved without the Regularity lemma). In order to obtain a
perfect $C_4$-packing, we have to incorporate the exceptional vertices.

To make it simpler to deal with divisibility issues later on, for every odd~$i$
we will now choose a set~$X_i$ of~7 vertices of~$G$ which we can put in any of~$V'_i$ and~$V'_{i+1}$
without destroying the superregularity of $(V'_i,V'_{i+1})_{G'}$. More precisely, (\ref{eq:minR})
implies that the vertices~$i$ and~$i+1$ of~$R$ have a common neighbour, $j$ say.
Recall that both $(V'_i,V'_j)_{G'}$ and $(V'_{i+1},V'_j)_{G'}$ are $2\eps$-regular and have
density at least $d/2$. So almost all vertices in~$V'_j$ have at least $(d/2-2\eps)m'$ neighbours in
both~$V'_i$ and~$V'_{i+1}$. Let~$X_i\subseteq V'_j$ be a set of~7 such vertices.
Clearly, we may choose the sets $X_i$ disjoint for distinct odd~$i$. Remove all the vertices
in $X_1\cup X_3\cup\dots\cup X_{k-1}=:X$ from the clusters they belong to. By removing at most
$|X|k\le 7k^2$ further vertices and adding them to the exceptional
set we may assume that the subclusters $V''_i\subseteq V'_i$ thus obtained satisfy
$|V''_1|=\dots =|V''_k|=:m''$. (The vertices in~$X$ are not added to~$V_0$.) Note that we now have
$$
|V_0|\le 4\eps n+7k^2\le 5\eps n.
$$
Consider any vertex $x\in V_0$. Call an odd~$i$ \emph{good for~$x$}
if~$x$ has at least~$\eta^2 m''$ neighbours in both~$V''_i$ and~$V''_{i+1}$ (in the graph~$G'$).
Then the number~$g_x$ of good indices satisfies
$$ (1/2+\eta/2)n\le d_{G'}(x)-|V_0|-|X|\le  2g_x m''+(k/2-g_x)(1+\eta^2) m''\le 2g_x m''+ (1+\eta^2)n/2,
$$
which shows that%
     \COMMENT{Here we need that $\eta<1/2$}
$g_x\ge \eta k/8=\eta |M|/4$. Since $|V_0|/(\sqrt{\eps} m'')\le \eta |M|/4$,
this implies that we can assign each $x\in V_0$ to an odd index~$i$
which is good for~$x$ in such a way that to each odd~$i$ we assign at most 
$\sqrt{\eps} m''$ exceptional vertices. Now consider any matching edge $i(i+1)\in M$.
Add each exceptional vertex assigned to~$i$ to~$V'_i$ or~$V'_{i+1}$ so that the
sizes of the sets $V^*_i\supseteq V''_i$ and $V^*_{i+1}\supseteq V''_{i+1}$ obtained in this
way differ by at most~1. It is easy to check that the bipartite subgraph $(V^*_i,V^*_{i+1})_{G'}$
of~$G'$ is still $(2\sqrt{\eps}, d/8)$-superregular.

Since the vertices in~$X_i$ can be added to
any of~$V^*_i$ and~$V^*_{i+1}$ without destroying the superregularity of $(V^*_i,V^*_{i+1})_{G'}$, we
could now apply the Blow-up lemma to find a $C_4$-packing of $G'[V^*_i\cup V^*_{i+1}\cup X_i]$
which covers all but at most 3 vertices (and so altogether these packings would form a
$C_4$-packing of $G$ covering all but at most $3k$ vertices of~$G$). To ensure the existence
of a perfect $C_4$-packing, we need to make $|V^*_i\cup V^*_{i+1}\cup X_i|$ divisible by~4 for every odd~$i$.
We will do this for every $i=1,3,\dots,k-1$ in turn by shifting the remainders $\mod 4$
along the path~$P$. More precisely, suppose that $|V^*_1\cup V^*_{2}\cup X_1|\equiv a\mod 4$
where $0\le a<4$. Choose $a$ disjoint copies of~$C_4$, each having~1 vertex in~$V^*_2$, 2 vertices in~$V^*_3$
and~1 vertex in~$V^*_4$. Remove the vertices in these copies from the clusters they belong to
and still denote the subclusters thus obtained by~$V^*_i$.
(Each such copy of~$C_4$ can be found greedily using that both $(V^*_2,V^*_3)_{G'}$ and $(V^*_3,V^*_4)_{G'}$
are still $2\sqrt{\eps}$-regular and have density at least $d/8$. Indeed, to find the first copy,
pick any vertex $x\in V^*_2$ having at least $(d/8-2\sqrt{\eps})|V^*_3|$ neighbours in~$V^*_3$.
The regularity of $(V^*_2,V^*_3)_{G'}$ implies that almost all vertices in~$V^*_2$ can play
the role of~$x$. The regularity of $(V^*_3,V^*_4)_{G'}$ now implies that its bipartite
subgraph induced by the neighbourhood of~$x$ in~$V^*_3$ and by~$V^*_4$ has density at
least $d/8-2\sqrt{\eps}$. So there are many vertices $y\in V^*_4$ which have at least~2 neighbours
in $N_{G'}(x)\cap V^*_3$. Then~$x$ and~$y$ together with~2 such neighbours form a copy of~$C_4$.)
Now $|V^*_1\cup V^*_2\cup X_1|$ is divisible by~4. Similarly, by removing at most~3 further copies of~$C_4$,
each having~1 vertex in~$V^*_4$, 2 vertices in~$V^*_5$ and~1 vertex in~$V^*_6$ we can
achieve that $|V^*_3\cup V^*_4\cup X_3|$ is divisible by~4.
Since $n=|G|$ is divisible by~$4$ we can continue in this way to achieve that  
$|V^*_i\cup V^*_{i+1}\cup X_i|$ divisible by~4 for every odd~$i$.

Recall that before we took out all these copies of~$C_4$, for every odd~$i$ the sizes of~$V^*_i$
and~$V^*_{i+1}$ differed by at most~1. Thus now these sizes differ (crudely) by at most~7.
But every vertex $x\in X_i$ can be added to both~$V^*_i$ and~$V^*_{i+1}$ without destroying
the superregularity. Add the vertices from~$X_i$ to~$V^*_i$ and~$V^*_{i+1}$
in such a way that the sets $V^\diamond_i\supseteq V^*_i$ and $V^\diamond_{i+1}\supseteq V^*_{i+1}$
thus obtained have equal size. (This size must be even since $|V^*_i\cup V^*_{i+1}\cup X_i|$
is divisible by~4.)
It is easy to check that $(V^\diamond_i,V^\diamond_{i+1})_{G'}$ is still
$(3\sqrt{\eps},d/9)$-superregular. Thus we can apply the Blow-up lemma (Lemma~\ref{blowup})
to obtain a perfect $C_4$-packing in $(V^\diamond_i,V^\diamond_{i+1})_{G'}$.
The union of all these packings (over all odd~$i$) together with the~$C_4$'s we have chosen
before form a perfect $C_4$-packing of~$G$.
\endproof
%The proof can be generalized to work for arbitrary graphs~$F$ instead of~$C_4$
%(where we replace 
%This gives a proof of the main result in~\cite{AY96} which is shorter than the original one

\section{Acknowledgment}
We would like to thank Demetres Christofides, Nikolaos Fountoulakis and Andrew Treglown for their
comments on an earlier version of this manuscript.

{\footnotesize
\bigskip\obeylines\parindent=0pt
Daniela K\"uhn \& Deryk Osthus
School of Mathematics
Birmingham University
Edgbaston
Birmingham B15 2TT
UK
{\it E-mail addresses}: {\tt \{kuehn,osthus\}@maths.bham.ac.uk}
}

\end{document}